\documentclass[11pt]{amsart}
\usepackage{tabularx,booktabs,tikz}
\usepackage{caption}
\usepackage{amsmath}
\usepackage{amsfonts}

\usepackage{amscd}
\usepackage{amsthm}
\usepackage{amssymb} \usepackage{latexsym}
\usepackage{eufrak}
\usepackage{euscript}
\usepackage{epsfig}
\usepackage{graphics}
\usepackage{array}
\usepackage{enumerate}
\usepackage{dsfont}
\usepackage{color}
\usepackage{wasysym}
\usepackage{hyperref}
\usepackage{pdfsync}

\newcommand{\bel}[1]{\begin{equation}\label{#1}}

\newcommand{\be}{\begin{equation}}

\newcommand{\ba}{\begin{eqnarray}}
\newcommand{\ea}{\end{eqnarray}}

\newcommand{\qe}{\end{equation}}
\newcommand{\rv}{{\mathbb R}^V}
\newcommand{\N}{{\mathbb N}}

\newcommand{\R}{{\mathbb R}}

\newcommand{\om}{\Omega}

\newcommand{\g}{\mathcal{G}}
\newcommand{\h}{\mathcal{H}}

\newcommand{\de}{{\Delta}}
\newcommand{\dO}{{\delta\Omega}}

\newcommand{\Hmm}[1]{\leavevmode{\marginpar{\tiny%
$\hbox to 0mm{\hspace*{-0.5mm}$\leftarrow$\hss}%
\vcenter{\vrule depth 0.1mm height 0.1mm width \the\marginparwidth}%
\hbox to
0mm{\hss$\rightarrow$\hspace*{-0.5mm}}$\\\relax\raggedright #1}}}

\newtheorem{theorem}{Theorem}[section]

\newtheorem{lemma}[theorem]{Lemma}
\newtheorem{corollary}[theorem]{Corollary}

\newtheorem{remark}[theorem]{Remark}

\newtheorem{prop}[theorem]{Proposition}
\newtheorem{problem}[theorem]{Problem}

\newtheorem{example}[theorem]{Example}

\newcommand{\tm}{\begin{theorem}}
\newcommand{\tmd}{\end{theorem}}
\newcommand{\co}{\begin{corollary}}
\newcommand{\cod}{\end{corollary}}
\newcommand{\prp}{\begin{prop}}
\newcommand{\prpd}{\end{prop}}
\newcommand{\pf}{\begin{proof}}
\newcommand{\pfd}{\end{proof}}
\newcommand{\rmk}{\begin{remark}}
\newcommand{\rmkd}{\end{remark}}
\newcommand{\ex}{\begin{example}}
\newcommand{\exd}{\end{example}}
\newcommand{\pr}{\begin{problem}}
\newcommand{\prd}{\end{problem}}

\begin{document}

\title[Steklov problems on trees]{Upper bounds for the Steklov eigenvalues on trees}

\author{Zunwu He}
\address{Zunwu He:School of Mathematical Sciences, Fudan University, Shanghai 200433, China}
\email{hzw@fudan.edu.cn}
%\address{Yohji Akama: Mathematical Institute, Graduate School of Science, Tohoku University,
%Sendai, 980-0845, Japan}
%\email{akama@math.tohoku.ac.jp}

\author{Bobo Hua}
\address{
Bobo Hua: School of Mathematical Sciences, LMNS, Fudan University, Shanghai 200433, China;  Shanghai Center for Mathematical Sciences, Jiangwan Campus, Fudan University, No. 2005 Songhu Road, Shanghai 200438, China.}
\email{bobohua@fudan.edu.cn}

%\author{Yanhui Su}
%\email{suyh@fzu.edu.cn}
%\address{Yanhui Su: College of Mathematics and Computer Science, Fuzhou University, Fuzhou 350116, China}

\begin{abstract}

%In this paper, we study the bounds for discrete Steklov eigenvalues on trees via geometric quantities. For a finite tree, we prove a sharp upper bound for the first nonzero Steklov eigenvalue by the reciprocal of the size of the boundary. This can be improved to the reciprocal of the size of the set of vertices for those trees with degree at least three for interior vertices.

%On the other hand, we obtain a sharp universal upper bound for the finite tree, in terms of the diameter of the tree.
%Moreover, we prove similar estimates for higher order Steklov eigenvalues.

In this paper, we study the upper bounds for discrete Steklov eigenvalues on trees via geometric quantities. For a finite tree, we prove sharp upper bounds for the first nonzero Steklov eigenvalue by the reciprocal of the size of the boundary and the diameter respectively.
We also prove similar estimates for higher order Steklov eigenvalues.
\end{abstract}
\maketitle

%\tableofcontents
%\tableofcontents
%Mathematics Subject Classification 2010: 05C10, 31C05.

\par
\maketitle

\bigskip

%\begin{enumerate}
%\item all first names are omitted.
%\end{enumerate}

\section{introduction}

The spectra of linear operators are important objects in Riemannian geometry, partial differential equations, graph theory, mathematical physics and so on. The Steklov problem is a classical eigenvalue problem exhibiting interesting
interactions between analysis and geometry; see e.g.  \cite{Stekloff02,Szeg54,Hersch75,Bandle80,Escobar97,Colbois11,Kuznetsov14,Jammes15,Girouard17}.
%\cite{Szeg\"{o}54},

Given a compact orientable Riemannian manifold $(M,g)$ with smooth boundary $\partial M$, the Steklov problem on $(M,g)$ reads as
\begin{align*}
\left\{
\begin{array}{lr}
\Delta f(x)=0,\ \ x\in M,&
\\\frac{\partial f}{\partial n}(x)=\lambda f(x),\ x\in \partial M,&
\end{array}
\right.
\end{align*}
where $\Delta$ is the Laplace-Beltrami operator on $(M,g)$ and $\frac{\partial}{\partial n}$ is the outward normal derivative along $\partial M$. The spectrum of {the} Steklov problem on $(M,g)$ coincides with that of the following Dirichlet-to-Neumann operator \cite{Kuznetsov14},
\begin{align*}
\Gamma:&H^{\frac{1}{2}}(\partial M)\longrightarrow H^{-\frac{1}{2}}(\partial M)
\\&f\longmapsto\Gamma f:=\frac{\partial \hat f}{\partial n},
\end{align*}
where $\hat f$ is the harmonic extension to $M$ of $f.$
%\cite{Ne\v{c}as12}
It is well-known that the Dirichlet-to-Neumann operator is a first-order elliptic pseudo-differential operator \cite{Taylor11}, which is self-adjoint and non-negative. Hence the spectrum is discrete, and can be ordered as $$0=\lambda_1<\lambda_2\leq \lambda_3\leq\cdots\nearrow\infty.$$ Here $\lambda_2$ is called the first (nonzero) Steklov eigenvalue.

In 1954, Weinstock \cite{Weinstock54} proved that the disk maximizes the first Steklov eigenvalue for simply connected planar domains of fixed perimeter.
For bounded Lipschitz domains of fixed volume in
$\R^n$, Brock \cite{Brock01} proved that the ball maximizes the the first  Steklov eigenvalue; see also \cite{Fraser11, Girouard12}. Escobar \cite{Escobar97,Escobar99,Escobar00}
systematically estimated the first Steklov eigenvalues using geometric quantities.
Colbois, El Soufi and Girouard \cite{Colbois11} proved an interesting result for the Steklov eigenvalues of a domain $\Omega$ in space forms, i.e. $\R^n,$ the hyperbolic space $\mathbb{H}^n$ and the sphere $\mathbb{S}^n,$
\begin{equation}\label{eq:h1}\lambda_k(\Omega)\leq \frac{C_n k^{\frac2n}}{\mathrm{Area}(\partial \Omega)^{\frac{1}{n-1}}}.\end{equation}

In order to detect spectral properties of Riemannian manifolds, Colbois et al. investigated the Steklov problem on some discretizations of manifolds \cite{Colbois18}. The second author, Huang and Wang \cite{Hua17}, and Hassannezhad and Miclo \cite{Miclo2017}, introduced the Cheeger-type isoperimetric constants to estimate the first Steklov eigenvalues on graphs independently. Perrin \cite{Perrin19} proved the lower bound estimate of the first Steklov eigenvalue for graphs; see \cite{ShiYu20a} for lower bounds using discrete curvatures. Han and the second author \cite{Han19} obtained the upper bound estimate of the first Steklov eigenvalue for finite subgraphs in integer lattices.  Recently, Perrin \cite{Perrin20} proved an upper bound estimate of the first Steklov eigenvalues for subgraphs in Cayley graphs of discrete groups of polynomial growth. Note that he essentially adopted the volume doubling property for these graphs, which fails for the trees of exponential growth.
See e.g. \cite{HuaHuangWang18,ShiYu19a,ShiYu19b,ShiYu20b} for other developments.

In this paper, we study eigenvalue estimates for the Steklov problem on trees. Trees can be regarded as discrete counterparts of Hadamard manifolds, i.e. simply-connected Riemannian manifolds of non-positive sectional curvature.

Let $\mathcal{G}=(V,E)$ be a finite tree. We denote by $\dO$ the set of boundary vertices, i.e.  those vertices of degree one, and by $\Omega=V\setminus \dO$ the interior vertices.
We write $|\cdot|$ the cardinality of a subset $(\cdot).$ We define the Steklov problem on
the pair $(\mathcal{G},\delta\Omega);$ see next section for definitions.
%\setminus
%In this paper, we focus mainly on the upper bounds of all Steklov eigenvalues for trees on account of two reasons. First, the lower bound of the first  Steklov eigenvalue for trees coincides with that for graphs and can be obtained for both cases. Second, Steklov eigenvalues have multiplicities.
For any function $0\neq f\in\rv,$ the  Rayleigh quotient of $f$ is defined as
\begin{align}\label{ray}
R(f):=\dfrac{\sum\limits_{(x,y)\in E}(f(x)-f(y))^2}{\sum\limits_{x\in \delta\Omega}f^2(x)},
\end{align}
where the right hand side is understood as $+\infty$ if $f|_{\dO}\equiv 0.$% $R(f)=\infty$ as $f(x)=0$ for any $x\in\dO$.

The variational characterization of the Steklov eigenvalues is given as
\begin{align}\label{mlk}
\lambda_k=\min_{\substack{W\subset \mathbb{R}^V,dimW=k-1\\W\perp 1_{\delta\Omega}}}\max\limits_{0\neq f\in W}
R(f),
\end{align}
where $1_{\delta\Omega}$ is the characterastic function on $\delta\Omega$, i.e. $f(x)=1$ if $x\in \delta\Omega$ and $f(x)=0$ for otherwise, and $W\perp 1_{\delta\Omega}:=\{g\in \rv: \sum_{x\in\delta \Omega} g(x)=0\}.$

The idea of getting the upper bound for the Escobar Cheeger-type constant in \cite{Hua17} motivates the authors to obtain the upper bounds for the first Steklov eigenvalue of finite trees, see Theorem \ref{lam2,1}, Theorem \ref{ulamk}.
A key observation is the following combinatorial property of trees: there remain two connected components if any edge of a tree is removed. We adopt it to divide the tree into two subgraphs with comparable sizes of boundaries. Then we construct a nonzero function $f$ orthogonal to $1_\dO,$ which is constant on each of the subgraphs with comparable constants. The first main result follows from the variational principle \eqref{mlk}.

 %such that $f$ is a constant on one subgraph, and another comparable constant on the other subgraph. By the formula (\ref{mlk}), we have the first main result.
\tm\label{lam2,1} Let $\mathcal{G}=(V,E)$ be a finite tree with the boundary $\delta\Omega$ such that the degree is bounded above by $D.$
Then we have
\begin{align}\label{l2,1}
\lambda_2\leq \dfrac{4(D-1)}{|\delta\Omega|}.
\end{align}
\tmd
%\rmk %\Hmm{I don't understand}
%The upper bound in (\ref{l2,1}) can be improved in terms of the vertices set $V$. The interested reader is referred to see
%\begin{enumerate}
%\item Note that our upper bound estimate is of order $\frac{1}{|\dO|},$ compared with the result \eqref{eq:h1} for domains in the hyperbolic space $\mathbb{H}^n$ by Colbois, El Soufi and Girouard, which is of order  $\left(\frac{1}{\mathrm{Area}(\partial \Omega)}\right)^{\frac{1}{n-1}}.$
We present a theorem to show that the estimate cannot be improved to $\frac{C(D)}{|V|}$ for general trees of bounded degree; see Theorem \ref{g3}.
%\end{enumerate}
%\rmkd

We say that a tree is of degree at least three if the degree of any interior vertex is at least three. Trees of degree at least three can be regarded as discrete counterparts of Hadamard manifolds with strictly negative curvature.
Next, we improve the above estimates by using the size of the set of vertices for such trees.

% for trees whose degree for any interior vertex is at least three, which

%However, the upper bound can be characterized in terms of $V$ under some mild constraint. It is formulated as follows.
\tm\label{uv}
For any finite tree $\mathcal{G}=(V,E)$ with the boundary $\delta\Omega$, and $\mathcal{G}$ has bounded degree $D$. If $\deg(x)\geq 3$ for any $x\in\om$, then we get
\begin{align}\label{lam2,2}
\lambda_2\leq \dfrac{8(D-1)}{|V|+2}.
\end{align}
\tmd

%Another description for the upper bound of the first  Steklov eigenvalue is the following.
Moreover, we estimate the first Steklov eigenvalue using the diameter of the tree.
\tm\label{2d}  For a finite tree $\mathcal{G}=(V,E)$ with the boundary $\delta\Omega$,
\begin{align}\label{lam2,3}
\lambda_2\leq \dfrac{2}{L},
\end{align}
where $L$ is the diameter of $\mathcal{G}$. Moreover, a necessary condition to attain the upper bound is that if $L$ is even, then the above $\g_k$ is exactly one point for any $1\leq k\leq L-1, k\neq \frac{L}{2}$ and if $L$ is odd, then $\g$ is exactly a path.
\tmd

%An interesting question on the upper bound is established.
%\pr\label{q1}
%\red{How to describe the finite trees that realize the upper bound in Theorem~\ref{2d}?}

%\prd

As a consequence, we prove the following result.
\co\label{2tend0}
For a sequence of finite trees $\mathcal{G}_n=(V_n,E_n)$ with uniformly bounded degree $D$. If $|V_n|\to \infty,$ $n\to\infty,$ then
$$\lim\limits_{n\to \infty}\lambda_2(\g_n)=0.$$
\cod

Aiming for estimating the higher order Steklov eigenvalues of trees, we further develop the previous ideas. First, we divide the tree into $k-1$ subgraphs with comparable sizes of boundaries. Second, we divide each subgraph into two subsets to construct a suitable function associated to each subgraph. These $k-1$ functions form a $(k-1)-$dimensional spaces $W$ in $\rv$. We estimate the higher order Steklov eigenvalue $\lambda_k$ using the Rayleigh quotients of functions in $W,$ which are properly controlled by our construction.

%By the construction, it can be showed that the maxiaml Rayleigh quotient of functions in $W$ can be attained by that of the previous $k-1$ functions. Third, the Rayleigh quotients of these $k-1$ functions are not hard to estimate. Thus we have
\tm\label{ulamk}
For any finite tree $\mathcal{G}=(V,E)$ with the boundary $\delta\Omega$, and $\mathcal{G}$ has bounded degree $D$. Then for any $3\leq k\leq |\dO|$, we have
\begin{align}\label{ulamk1}
\lambda_k\leq \dfrac{8(D-1)^2(k-1)}{|\delta\Omega|}.
\end{align}
\tmd

Note that the constants in the above estimate involve only $k$ and $D,$ which is linear in $k$, compared with the result \eqref{eq:h1} for the hyperbolic space $\mathbb{H}^n,$ which is of order $k^\frac{2}{n}.$

The estimate can be improved using the size of the set of vertices for trees of degree at least three.
%We also have the following result which relates the upper bound of $k$-th Steklov eigenvalue to the vertices set.
\co\label{kV}
Let $\mathcal{G}=(V,E)$ be a finite tree with the boundary $\delta\Omega$ such that the degree is bounded above by $D.$ If $\deg(x)\geq 3$ for any $x\in\om$, then
\begin{align*}
\lambda_k\leq \dfrac{16(D-1)^2(k-1)}{|V|+2}.
\end{align*}
\cod

\textbf{Acknowledgements.} B.H. is supported by NSFC, grants no.11831004 and no. 11926313. The authors very appreciate anonymous referees for providing many
valuable suggestions to improve the writing of the paper. The authors are also grateful to Prof. Thomas Schick and Prof. Chengjie Yu for many helpful discussions.

%\textbf{Data Available Statement.} The data used to support the findings of this study are available from the corresponding author upon request.

%\textbf{Data availability statements.} The data that support the findings of this study are available from the corresponding author upon reasonable request.

\section{Preliminaries}

For a simple, undirected graph $\mathcal{G}=(V,E),$ two vertices $x,y$ are called neighbors, denoted by $x\sim y,$ or $(x,y)\in E$ if there is an edge connecting $x$ and $y.$ We say $\mathcal{G}$ is connected if for any $x,y\in V$, there is a path $x=x_0\sim x_1\sim \cdots \sim x_n=y$ connecting $x$ and $y$ for some $n\in\N.$ For any $x\in V,$ we denote by $\deg(x)$ the vertex degree of $x.$ If $\deg(x)\leq D$ for any $x\in V$ and some positive integer $D$, we say the graph has bounded degree $D$. The combinatorial distance $d$ on the graph is defined as, for any $x,y\in V$ and $x\neq y,$ $$d(x,y):=\inf\{n\in\N: \exists\{x_{i}\}_{i=1}^{n-1}\subset V, x\sim x_1\sim\cdots \sim x_{n-1}\sim y \}.$$ % We denote by $|\cdot |$ the counting measure on $V,$ i.e. for any $V_1\subset V,$ $|V_1|$ denotes the cardinality of vertices in $V_1$.

%We introduce the boundary and the interior of $\mathcal{G}$ in the following way. The \emph{boundary} $\delta\Omega$ is the set of vertices with degree one, and the \emph{interior} $$\Omega:=V\backslash\delta\Omega.$$
From now on, we always consider finite trees $\mathcal{G}$ of bounded degree $D$ with at least three vertices to exclude that $\Omega=\emptyset$. Denote by $(\mathcal{G},\delta\Omega)$ the graph $\mathcal{G}=(V,E)$ with the boundary $\delta\Omega$ and $\Omega=V\setminus\dO.$ Recall that $\delta\Omega$ is the set of vertices of degree one in $\g$.

For a subtree $\g_1=(V_1,E_1)$ in $\g$, we define the relative boundary to be the subset of $\dO$ in $V_1$, denoted by $\dO\cap\g_1$ or $\dO(\g_1)$.

Note that there is no edge connecting two boundary vertices and $\Omega$ is connected.  It is clear that this graph can be embedded in a homogeneous tree with degree no less than $D$.

Denote by $\mathbb{R}^{V}$ the vector space of all real functions on $V$ over $\mathbb{R}$, and it can be equipped with $\ell^2$-inner product $(\cdot,\cdot)$. For any $f,g\in \mathbb{R}^{V}$, one can set $(f,g):=\sum\limits_{x\in V}f(x)g(x)$, then $(\mathbb{R}^{V},(\cdot,\cdot))$ is a Hilbert space. For $f\in \mathbb{R}^V$, one can define a Laplace operate $\Delta$ on $\mathbb{R}^{V}$ such that
\begin{align}\label{laplace}
(\Delta f)(x):=\sum\limits_{y\sim x}(f(x)-f(y))
\end{align}
A function $g\in \mathbb{R}^V$ is called harmonic if $(\Delta g)(x)=0$ for any $x\in \Omega$.

In analogy to the Riemannian setting, with a fixed orientation for edges, we introduce the gradient of a function $f\in\rv$
\begin{align*}
& \nabla f:\ E\longrightarrow \R
\\&e=(x,y)\longmapsto \nabla_e f:=f(y)-f(x).
\end{align*}

One can define the outward normal derivative operator
\begin{align*}
\frac{\partial}{\partial n}:&\mathbb{R}^V\longrightarrow \mathbb{R}^{\delta\Omega}
\\&f\longmapsto\frac{\partial f}{\partial n},
\end{align*}
where $\frac{\partial f}{\partial n}(x)=\sum\limits_{y\in \Omega,y\sim x}(f(x)-f(y))$ for any $x\in \delta \Omega.$ Note that $\frac{\partial f}{\partial n}(x)=(\Delta f)(x)$, since there is no edge connecting two boundary vertices.

We introduce the Steklov problem on the pair $(\mathcal{G},\delta\Omega)$. For some nonzero function $f\in \mathbb{R}^V$ and some $\lambda\in \mathbb{R}$, the following equation holds.
\begin{equation}\label{Steklov 2}
\left\{
\begin{array}{lc}
\Delta f(x)=0, x\in \Omega&\\
\dfrac{\partial f}{\partial n}(x)=\lambda f(x), x\in \delta\Omega.
\end{array}
\right.
\end{equation}
The above $\lambda$ is called the Steklov eigenvalue of the graph with boundary $(\mathcal{G},\delta\Omega)$, and $f$ is the Steklov eigenfunction associated to $\lambda$.

In analogy to the Riemannian case, one can define the Dirichlet-to-Neumann operator in the discrete setting to be
\begin{align}
\Lambda:&\mathbb{R}^{\delta\Omega}\longrightarrow \mathbb{R}^{\delta\Omega}\nonumber
\\&f\longmapsto \Lambda f:=\frac{\partial \hat f}{\partial n},\label{dtn}
\end{align}
where $\hat f$ is the harmonic extension of $f$. Namely, $\hat f$ satisfies
\begin{equation}\label{dtn2}
\left\{
\begin{array}{lr}
\Delta \hat f(x)=0, x\in \Omega & \\
\hat f(x)=f(x), x\in \delta\Omega &
\end{array}
\right.
\end{equation}
It is well known that the Steklov eigenvalues in (\ref{Steklov 2}) are exactly the Dirichlet-to-Neumann eigenvalues in (\ref{dtn}), and the Steklov eigenfunctions in (\ref{Steklov 2}) are  the harmonic extensions of the corresponding eigenfunctions in (\ref{dtn}). The Dirichlet-to-Neumann operator is non-negative, self-adjoint. Since $V$ is finite, there are $|\dO|$ Steklov eigenvalues. We may arrange the Steklov eigenvalues in the following way:
\begin{align}\label{steklov eigen}
0=\lambda_1\leq \lambda_2\leq \cdots\leq \lambda_{|\dO|}.
\end{align}

In fact, $\Omega$ is connected, then $\lambda_2>0$ and $\lambda_{|\dO|}\leq 1$. Moreover, $(f_i,f_j)=0$ if $\lambda_i\neq \lambda_j$, where $f_i$ and $f_j$ are eigenfunctions with eigenvalues $\lambda_i,\lambda_j$ respectively. $\mathbb{R}^{\delta\Omega}$ consists of all the Dirichlet-to-Neumann eigenfunctions. Note that constant functions are the eigenfunctions of $\lambda_1=0$.

%For $f\in $
For any $0\neq f\in \rv$, recall that we define the Rayleigh quotient
\begin{align}\label{ray}
R(f)=\dfrac{\sum\limits_{(x,y)\in E}(f(x)-f(y))^2}{\sum\limits_{x\in \delta\Omega}f^2(x)}.
\end{align}
Some variational descriptions for $\lambda_k$ are
\begin{align}\label{lamk}
&\lambda_k=\min\limits_{W\subset \mathbb{R}^V,\ dimW=k}\max_{0\neq f\in W}R(f)
\\&\lambda_k=\min_{\substack{W\subset \mathbb{R}^V,dimW=k-1\\W\perp 1_{\delta\Omega}}}\max\limits_{0\neq f\in W}
R(f)\label{var k},
\end{align}
where $1_{\delta\Omega}$ is the characteristic function on $\delta\Omega$, i.e. $f(x)=1$ if $x\in \delta\Omega$ and $f(x)=0$ for otherwise.

\section{Upper bounds for the first  Steklov eigenvalue}
In the beginning of this section, we introduce a lemma which describes the relation between the boundary $\dO$ and the set of vertices $V$ of a finite tree $\g$. The simple proof of Lemma \ref{dv} is suggested by the referees.
\begin{lemma}\label{dv}
Let $\g=(V,E)$ be a finite tree with the boundary $\dO$ such that $\Omega\neq \emptyset.$ If $deg(x)\geq 3$ for any $x\in \Omega$, then the following holds
$$\frac{1}{2}|V|+1\leq |\dO|\leq |V|.$$
%Moreover, the upper and lower bounds for $|\dO|$ are sharp.
\end{lemma}

\pf
Since the tree is contractible and $deg(x)\geq 3$ for any $x\in \Omega$, we have
\begin{align}
&2(|V|-1)\geq 2|E|=\sum\limits_{x\in V}deg(x)
\\&\geq|\delta\Omega|+3|\Omega|
\\&=|\delta\Omega|+3(|V|-|\dO|)
\\&=3|V|-2|\dO|
\end{align}

This gives the proof.
\pfd

%The first main result is the upper bound of the first  Steklov eigenvalue, i.e. $\lambda_2$ in (\ref{steklov eigen}).
%\begin{theorem}\label{lam2,1}
%For any finite tree $\mathcal{G}=(V,E)$ with the boundary $\delta\Omega$, and $\mathcal{G}$ has bounded degree $D$. Then there exists a constant $c_1(D)$ only depending on $D$ such that
%\begin{align}\label{lam2}
%\lambda_2\leq \dfrac{c_1(D)}{|\delta\Omega|}.
%\end{align}
%\end{theorem}
Next, we give a combinatorial lemma which is key to approach the upper bound of the first Steklov eigenvalue of finite trees.
\begin{lemma}\label{2tree}
Let $\g=(V,E)$ be a finite tree of bounded degree $D$ with boundary $\dO$, then there exists a subtree $\h$ in $\g$ such that
$$\dfrac{|\h\cap \dO|}{|\dO|}\geq \dfrac{1}{2(D-1)},$$
 where $\h$ is one of two subtrees obtained by removing one edge of $\g$.
\end{lemma}

\pf
Choose an arbitrary edge $e_1\in E$, and there are two connected components $\Gamma_1,\Gamma_2$ if we remove $e_1$ in $\mathcal{G}$. If $\dfrac{|\Gamma_1\cap \delta\Omega|}{|\delta\Omega|}\leq \dfrac{1}{2}$, $\dfrac{|\Gamma_2\cap \delta\Omega|}{|\delta\Omega|}\leq \dfrac{1}{2}$, then $\dfrac{|\Gamma_1\cap \delta\Omega|}{|\delta\Omega|}= \dfrac{|\Gamma_1\cap \delta\Omega|}{|\delta\Omega|}=\dfrac{1}{2}$, we finish the proof.

Otherwise, we may assume $\dfrac{|\Gamma_2\cap \delta\Omega|}{|\delta\Omega|}>\dfrac{1}{2}$ and there are at most $D-1$ edges in $\Gamma_2$ adjacent to $e_1$. There are at most $D-1$ connected components if one removes these edges in $\Gamma_2$. If the cardinality of the relative boundary vertices of none of connected components is more than one half of $|\delta\Omega|$, we may assume an edge $e_2$ in $\Gamma_2$ is adjacent to $e_1$ and $\Gamma_3$ is the connected component resulted from removing $e_2$ in $\Gamma_2$, such that the cardinality of the relative boundary vertices of $\Gamma_3$ is maximal among those of above connected components in $\Gamma_2$. Then we have

\begin{align*}
&\dfrac{|\Gamma_3\cap \delta\Omega|}{|\delta\Omega|}\leq \dfrac{1}{2},
\end{align*}
and
\begin{align*}
\dfrac{|\Gamma_3\cap \delta\Omega|}{|\delta\Omega|}
=\dfrac{|\Gamma_3\cap \delta\Omega|}{|\Gamma_2\cap \delta\Omega|}\dfrac{|\Gamma_2\cap \delta\Omega|}{|\delta\Omega|}
\geq \dfrac{1}{D-1}\cdot \dfrac{1}{2}.
\end{align*}

Otherwise, the cardinality of the relative boundary vertices of some connected component is more than one half of $|\delta\Omega|$. For simplicity, we may denote by $\Gamma_3$ this connected component, $e_2$ the corresponding edge in $\Gamma_2$. And there are at most $D-1$ edges in $\Gamma_3$ adjacent to $e_2$. There are at most $D-1$ connected components if one removes these edges in $\Gamma_2$.

Repeating the above procedures and since $\g$ is finite, we can get a finite sequences of subgraphs $\Gamma_2\supseteq \Gamma_3\supseteq \cdots \supseteq \Gamma_i$ and corresponding consecutive edges $e_1,e_2,\cdots,e_{i-1}$ such that
\begin{align}
\dfrac{|\Gamma_i\cap \delta\Omega|}{|\delta\Omega|}>\dfrac{1}{2},\label{4,1}
\end{align}
and the cardinality of the relative boundary vertices of none of connected components, resulted from removing at most $D-1$ edges adjacent to $e_{i-1}$, is more than one half of $|\delta\Omega|$.

We may assume an edge $e_i$ in $\Gamma_i$ is adjacent to $e_{i-1}$ and $\Gamma_{i+1}$ is the connected component resulted from removing $e_i$ in $\Gamma_i$, such that the cardinality of the relative boundary vertices of $\Gamma_{i+1}$ is maximal among those of above connected components in $\Gamma_i$. Then using (\ref{4,1}) we have

\begin{align}\label{5,1}
&\dfrac{|\Gamma_{i+1}\cap \delta\Omega|}{|\delta\Omega|}\leq \dfrac{1}{2},
\end{align}
and
\begin{align}\label{5,4}
\dfrac{|\Gamma_{i+1}\cap \delta\Omega|}{|\delta\Omega|}
=\dfrac{|\Gamma_{i+1}\cap \delta\Omega|}{|\Gamma_i\cap \delta\Omega|}\dfrac{|\Gamma_i\cap \delta\Omega|}{|\delta\Omega|}
\geq \dfrac{1}{D-1}\cdot \dfrac{1}{2}.
\end{align}

Thus we get the conclusion by setting $\h:=\Gamma_{i+1}$.
\pfd

We are ready to prove the first main result of the paper, Theorem~\ref{lam2,1}.
\begin{proof}[Proof of Theorem~\ref{lam2,1}]
By (\ref{var k}) for $k=2$, we have $\lambda_2=\min\limits_{f\neq 0,\sum\limits_{x\in \delta\Omega}f(x)=0} R(f)$. Hence, one easily sees that $\lambda_2\leq R(f)$ for any $0\neq f\in \mathbb{R}^V$ with $\sum\limits_{x\in \delta\Omega}f(x)=0$.

%Choose an arbitrary edge $e_1\in E$, and there are two connected components $\Gamma_1,\Gamma_2$ if we remove $e_1$ in $\mathcal{G}$. If $\dfrac{|\Gamma_1\cap \delta\Omega|}{|\delta\Omega|}\leq \dfrac{1}{2}$, $\dfrac{|\Gamma_2\cap \delta\Omega|}{|\delta\Omega|}\leq \dfrac{1}{2}$, then $\dfrac{|\Gamma_1\cap \delta\Omega|}{|\delta\Omega|}= \dfrac{|\Gamma_1\cap \delta\Omega|}{|\delta\Omega|}=\dfrac{1}{2}$. We can assume $f\in \mathbb{R}^V$ satisfying $f(x)=1$ for $x\in \Gamma_1$, $f(x)=-1$ for $x\in \Gamma_2$, then
%\begin{align}\label{2,1}
%\lambda_2\leq R(f)=\dfrac{4}{|\delta\Omega|}.
%\end{align}

Let $\h$ be as in Lemma \ref{2tree}.

We set $f(x)=1-\dfrac{|\h\cap \delta\Omega|}{|\delta\Omega|}$ for $x\in \h$ and $f(x)=-\dfrac{|\h\cap \delta\Omega|}{|\delta\Omega|}$ for $x\notin \h$. It is clear that $(f,1_{\delta\Omega})=0$ and by (\ref{5,1}) and (\ref{5,4}), we get
\begin{align*}
\lambda_2&\leq R(f)
\\&\leq \dfrac{1}{((1-\frac{1}{2(D-1)})^2\frac{1}{2(D-1)}+(\frac{1}{2(D-1)})^2(1-\frac{1}{2(D-1)}))|\delta\Omega|}
\\&\leq \dfrac{4(D-1)}{|\delta\Omega|}.
\end{align*}

We complete the proof.
\end{proof}

\begin{remark}
The upper bound in (\ref{l2,1}) is sharp up to a universal constant, such as Example \ref{e1}.

\begin{example}\label{e1}
Let $T_D=(V, E)$ be a homogeneous tree of degree $D$ with the combinatorial metric. We assume that $\mathcal{G}_r=(V_r,E_r)$ is the induced subgraph on a ball $B(o,r)$ in $T_D$ of radius $r$, centralized at $o\in V_r$, then the boundary of $\mathcal{G}_r$ is $\delta\Omega_r=S(0,r)$, namely the sphere of radius $r$. The interior of $\mathcal{G}$ is $\Omega_r=B(o,r-1)$. See Figure~\ref{fig0}

\begin{figure}[htbp]
 \begin{center}
       \includegraphics[width=0.4\linewidth]{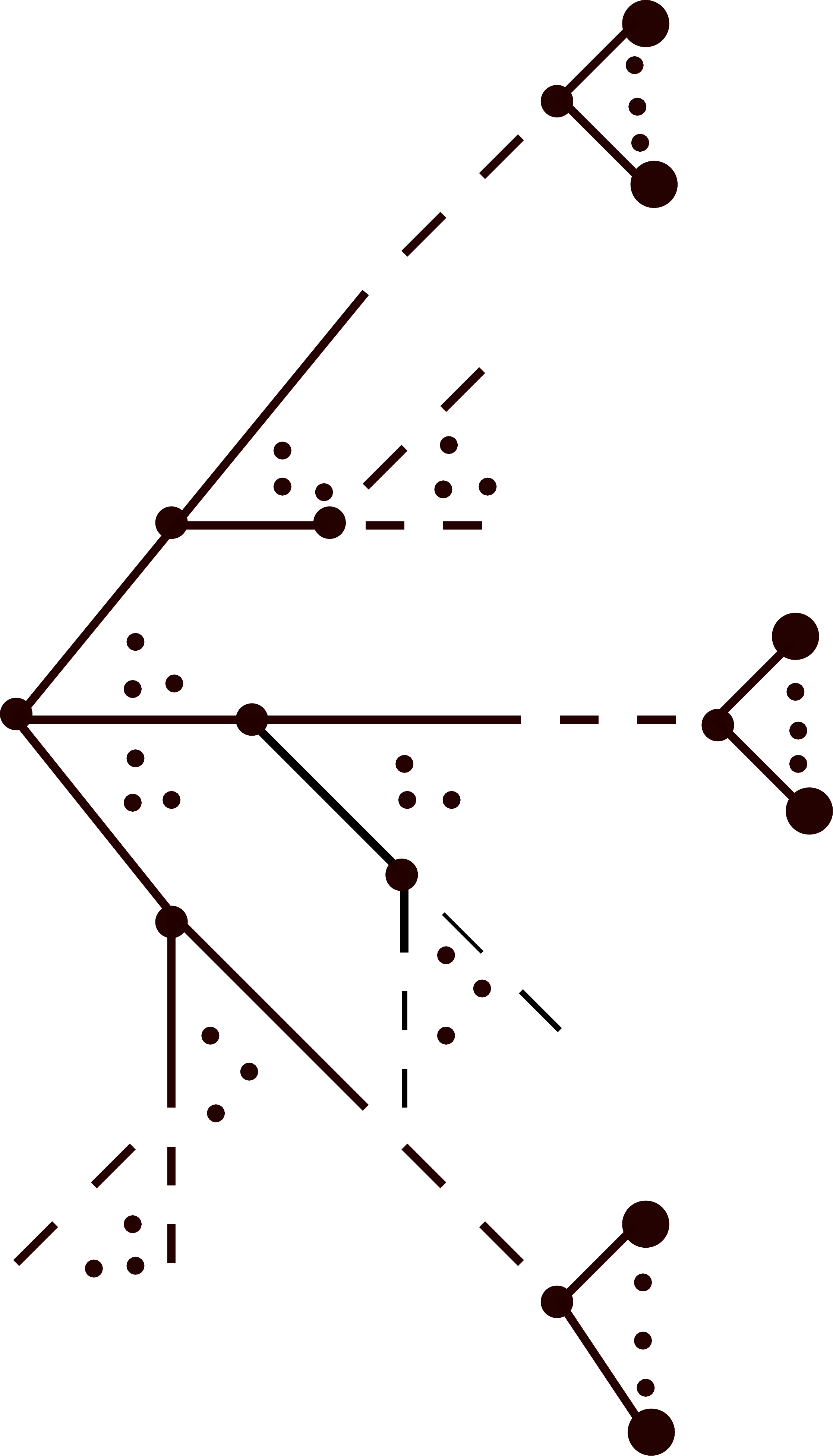}
  \caption{}
  %\small
 %This is a tree with $n$ vertices, $2$ boundary vertices. }
 \label{fig0}
 \end{center}

\end{figure}

Note that $|\delta\Omega_r|=D(D-1)^{r-1}$, $|\delta\Omega_r|\leq |V_r|\leq 2|\delta\Omega_r|$. A straightforward computation shows that $\lambda_2(\mathcal{G}_r)=\frac{1}{\sum\limits_{k=0}^{r-1}(D-1)^k}$, so we have
\begin{align*}
\dfrac{D}{2|\delta\Omega_r|}\leq \lambda_2(\mathcal{G}_r)\leq \dfrac{D}{|\delta\Omega_r|}.
\end{align*}

\end{example}
\end{remark}

Now we can prove the following result.
\pf[Proof of Theorem~\ref{uv}]
It is clear from Lemma \ref{dv} and Theorem \ref{lam2,1}.
\pfd

Concerning with Theorem \ref{lam2,1} and Theorem \ref{uv}, one may want to strengthen the upper bound in Theorem \ref{lam2,1} or remove the assumption on the degree of interior points in Theorem \ref{uv}. %It is stated as a problem.
\pr\label{p1}
Can (\ref{l2,1}) in Theorem~\ref{lam2,1} be improved to $\lambda_2\leq \dfrac{c(D)}{|V|}$ or the assumption $deg(x)\geq 3$ for all $x\in\Omega$ in Theorem \ref{uv} is necessary, where $c(D)>0$ is a constant depending only on $D$?
\prd

Unfortunately, the answer is that the improvement fails and the assumption is necessary in Problem \ref{p1} in general ; see the following theorem.
\tm\label{g3}
There exist a family of finite trees to show that the improvement fails and the assumption is necessary in Problem \ref{p1}. To be precise, let $\g_1=(V_1,E_1)$ be a subgraph generated by $B(o,l)$ in the homogeneous tree $T_3$ of degree three with the combinatorial metric $d$, $o\in V_1$ and $l\in \N$. We refine $\g_1$ to get $\g=(V,E)$ by adding $k$ vertices to the edge $(x,y)$, with $d(o,x)=k,d(o,y)=k+1$ and $1\leq k\leq l-1$.

\begin{figure}[htbp]
 \begin{center}
       \includegraphics[width=0.4\linewidth]{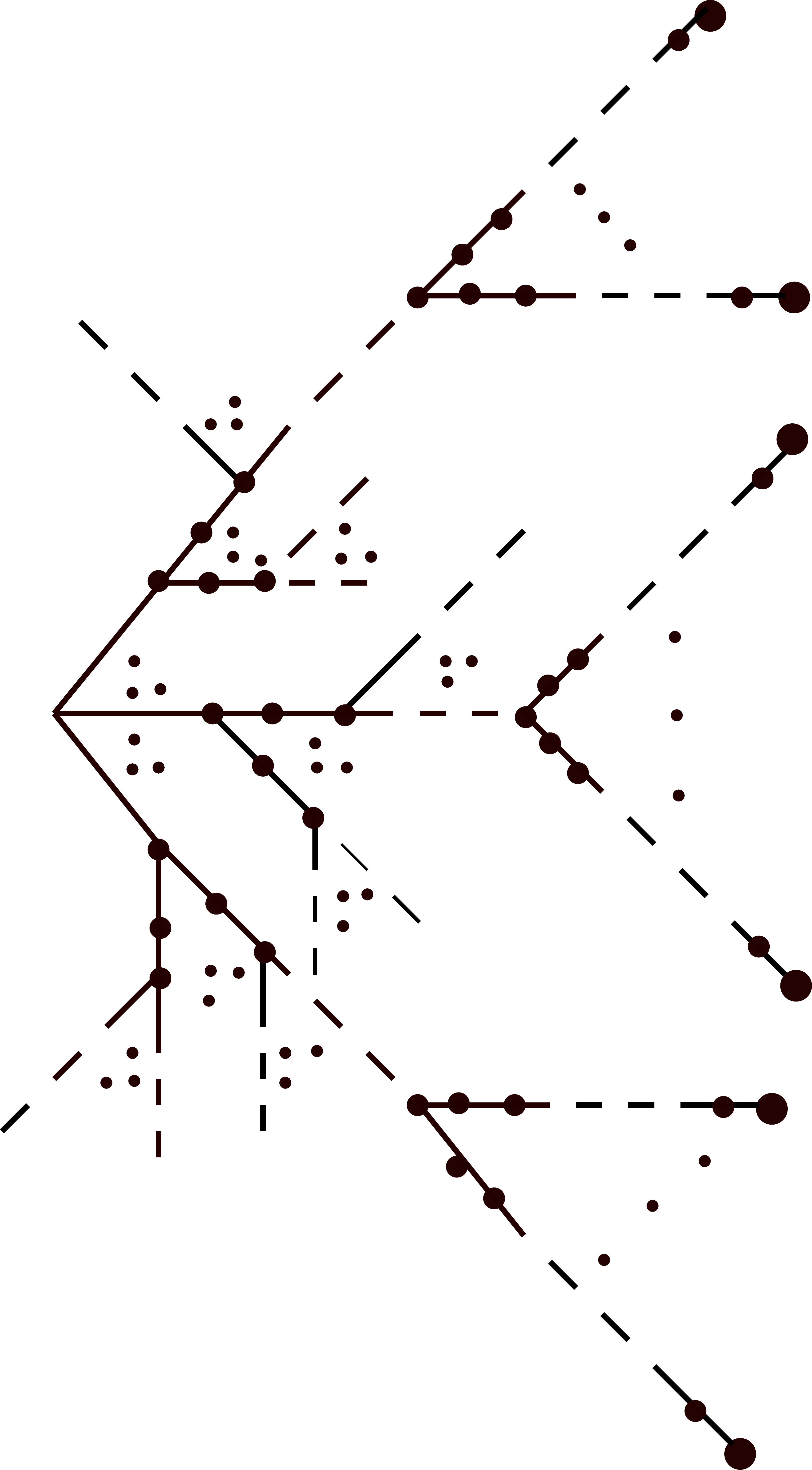}
  \caption{ }
 \label{fig3}
 \end{center}

\end{figure}

Denote by $\dO$, $\dO_1$ the boundary of $\g$, $\g_1$ respectively. See Figure~\ref{fig3}

Then we have
\begin{align*}
&|\dO|=|\dO_1|=3\cdot2^{l-1},
\\l|\dO|\leq&|V|=|V_1|+\sum\limits_{k=1}^{l-1}3\cdot2^k\cdot k\leq 2l|\dO|.
\end{align*}

Furthermore, it holds that
\begin{align}\label{g3,1}
\dfrac{3}{4|\dO|}\leq \lambda_2(\g)=\lambda_2\leq \dfrac{3}{2|\dO|}.
\end{align}
\tmd
\pf%[Proof of the inequality (\ref{g3,1}) in Example \ref{g3}]
We only show the inequality (\ref{g3,1}) in Theorem \ref{g3}. Denote by $f\in \rv$ one of the eigenfunctions of $\lambda_2$. There is a connected component $\g_2=(V_2,E_2)$ which does not contain $o$, if one removes an edge containing $o$. Similar to the case of $\g_1$, we may assume $f(z)\neq 0$ for any $z\in V_2$, $f(x)=f(y)$ for any two boundary vertex $x,y$ in $\g_2$ and $f(o)=0$. Since $\de f(x)=0$ for $x\notin \dO$ and $\de f(x)=\lambda_2f(x)$ for $x\in \dO,$ we get
\begin{align*}
\dfrac{3}{4|\dO|}=\dfrac{1}{2^{l+1}}\leq \lambda_2=\dfrac{1}{\sum\limits_{k=0}^{l}(l-k)2^k}=\dfrac{1}{2^{l+1}-l-2}\leq \dfrac{1}{2^l}=\dfrac{3}{2|\dO|}.
\end{align*}

\pfd

On the other hand, there exist some family of finite trees $\mathcal{G}_n$ with boundaries $|\delta\Omega_n|=2$ and $\lim\limits_{n\rightarrow \infty}\lambda_2(\mathcal{G}_n)=0$. Next, we prove another upper bound for $\lambda_2$ in terms of the diameter of the graph.

\begin{proof}[Proof of Theorem~\ref{2d}]
Using (\ref{lamk}) for $k=2$, we have
\begin{align}\label{22}
\lambda_2\leq \frac{\sum\limits_{(x,y)\in E}(f(x)-f(y))^2}{\sum\limits_{x\in \delta\Omega}f^2(x)}=R(f),
\end{align}
for any $0\neq f\in \mathbb{R}^V$ with $\sum\limits_{x\in \delta\Omega}f(x)=0$.

One easily shows that there exist two vertices $z,y\in \dO$ with $d(z,y)=L$, where $d$ is the combinatorial metric on $\g$. We may assume $x_0=z\sim x_1\sim x_2\cdots\sim x_L=y$.
For $x_k$ with $1\leq k\leq L-1$, there are three connected components if one removes edges $(x_{k-1},x_k),(x_k,x_{k+1})$. Denote by $\g_k=(V_k,E_k)$ the component which does not contain $x_0$ or $x_L$. Recall that the relative boundary $\dO(\Omega_k):=V_k\cap\dO$ and set $n_k:=|\dO(\Omega_k)|\in\N$.

\textbf{Claim:} One can choose $f\in \rv$ and $a_0,a_k\in\R$ such that
\begin{align}
f(x)=f(x_k)=a_0-k\cdot\dfrac{2a_0+\sum\limits_{i=1}^{L-1}n_ia_i}{L}\label{2d,1}
\end{align}
for $x\in V_k$ and $1\leq k\leq L-1$, $f(x_0)=a_0$ and $f(x_L)=-a_0-\sum\limits_{k=1}^{L-1}n_ka_k$.

\pf[Proof the claim]
Since $L$ is the diameter of $\g$, we have that $V=\{x_0,x_L\}\cup\bigcup\limits_{k=1}^{L-1}V_k$. On the other hand, we have the following equations.
\begin{equation*}
\left\{
\begin{array}{lr}
(L-2)a_0-(n_1+L)a_1-n_2a_2-n_3a_3-\cdots-n_{L-1}a_{L-1}=0 &\\
\cdots    &\\
(L-2k)a_0-n_1a_1-\cdots-n_{k-1}a_{k-1}-(n_k+L)a_k-n_{k+1}a_{k+1}-\cdots-n_{L-1}a_{L-1}=0 &\\
\cdots &\\
(L-2(L-1))a_0-n_1a_1-n_2a_2-\cdots-n_{L-2}a_{L-2}-(n_{L-1}+L)a_{L-1}=0&
\end{array}
\right.
\end{equation*}

Note that the above system of equations is $L-1$ linear equations with respect to $L$ variables $a_0,a_1,\cdots,a_{L-1}$. By the linear algebra theory, there exists a nonzero solution $(a_0,a-1,\cdots,a_{L-1})\neq 0$ satisfying the above system of equations.

Thus we conclude the claim.
\pfd

Moreover, we get that
$$f(x_{k-1})-f(x_{k})=f(x_k)-f(x_{k+1})=\dfrac{a_0+a_0+\sum\limits_{k=1}^{L-1}n_ka_k}{L}.$$

It is obvious that $f$ satisfies the condition for (\ref{22}). Then we deduce
\begin{align*}
\lambda_2\leq R(f)&=\dfrac{({\dfrac{a_0+a_0+\sum\limits_{k=1}^{L-1}n_ka_k}{L}})^2L}{a_0^2+(a_0+\sum\limits_{k=1}^{L-1}n_ka_k)^2+\sum\limits_{k=1}^{L-1}n_ka_k^2}\\
&\leq\frac{1}{L}\dfrac{{(a_0+a_0+\sum\limits_{k=1}^{L-1}n_ka_k)^2}}{a_0^2+(a_0+\sum\limits_{k=1}^{L-1}n_ka_k)^2}\\
&\leq\frac{2}{L}.
\end{align*}
The necessary condition to attain the upper bound follows from the above discussion.
\end{proof}

\rmk
The above estimate is universal which is independent of the upper bound of the degree information. The estimate is sharp; see trees attaining the bound in Example~\ref{exp:1}.
\rmkd

\rmk\label{line tree}
%From the above proof, one deduces that a necessary condition to attain the upper bound in Theorem~\ref{2d} is that if $L$ is even, then the above $\g_k$ is exactly one point
%for any $1\leq k\leq L-1, k\neq \frac{L}{2}$ and if $L$ is odd, then $\g$ is exactly a path.
The equivalent description of the rigidity of the isodiametric estimate is referred to our another paper for interested readers; see \cite{HH21}.
\rmkd

Now we give some examples attaining the upper bound \eqref{lam2,3} in Theorem~\ref{2d}. One can verify that via direct computation.%The verification is not hard and left to interested readers.

\ex\label{exp:1}
Let $L$ be an even number.  In the following examples, $\lambda_2=\frac{2}{L}.$

\begin{figure}[htbp]
 \begin{center}
       \includegraphics[width=0.6\linewidth]{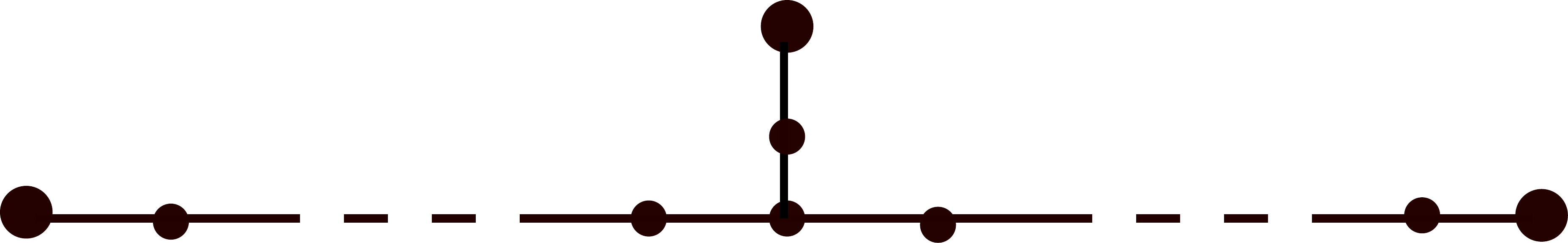}
  \caption{}
 \label{fig4}
 \end{center}

\end{figure}

Figure~\ref{fig4} is a finite tree with diameter $L$, $|\g_k|=1$ and $|\g_{\frac{L}{2}}|=3,$   where $1\leq k\leq L-1, k\neq \frac{L}{2},$ and $\frac{L}{2}\geq 2$.

\begin{figure}[htbp]
 \begin{center}
       \includegraphics[width=0.6\linewidth]{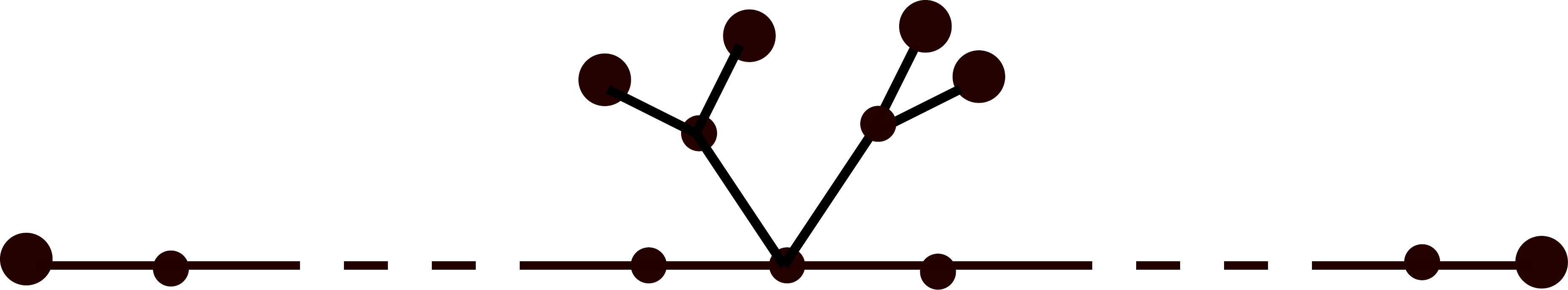}
  \caption{}
 %This is a tree with $n$ vertices, $2$ boundary vertices. }
 \label{fig5}
 \end{center}

\end{figure}

Figure~\ref{fig5} is a finite tree with diameter $L$, $|\g_k|=1$ and $|\g_{\frac{L}{2}}|=7$  , where $1\leq k\leq L-1,  k\neq \frac{L}{2},$ and $\frac{L}{2}\geq 3$.

\begin{figure}[htbp]
 \begin{center}
       \includegraphics[width=0.6\linewidth]{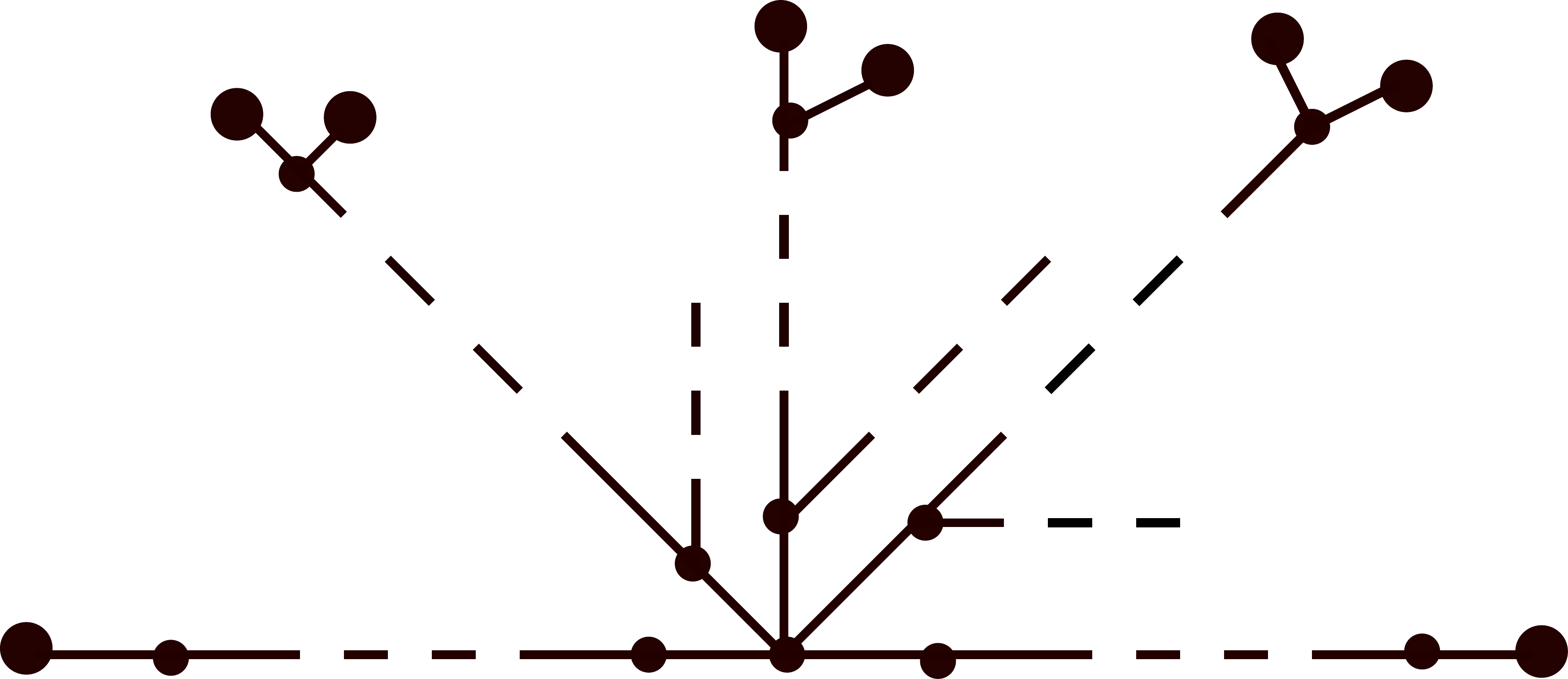}
  \caption{}
 %This is a tree with $n$ vertices, $2$ boundary vertices. }
 \label{fig6}
 \end{center}

\end{figure}

Figure~\ref{fig6} is a finite tree with diameter $L$, $|\g_k|=1$ and $\g_{\frac{L}{2}}$ is generated by the ball $B_{\hat L}$ of radius $\hat L$ in the regular tree of degree three, where $1\leq k\leq L-1, k\neq \frac{L}{2},$ and $\frac{L}{2}\geq \sum\limits_{i=0}^{\hat L-1}2^i$.

\exd

From the above examples, it seems difficult to characterize the equality case for \eqref{lam2,3} in Theorem~\ref{2d}.

%The answer of Question \ref{q1} seems quite complicated from the above examples. But we think that if a finite tree satisfying the necessary condition in Remark~\ref{line tree} has a sufficiently large (relative to the cardinality of the vertex set of $\g_{\frac{L}{2}}$) diameter $L$, then it realizes the upper bound.}

\prp\label{L}
Let $\g=(V,E)$ and $L$ be as above, then $L\geq 2\log_{D}\dfrac{|V|}{4}$.

\prpd

\pf
We keep the notations in the above proof.

If $L=2k$ for some positive integer $k$, then $d(x_0,x_k)=d(x_k,x_L)=k$. It is not hard to see that $|V|$ is maximal if $\g$ is a subgraph in $T_D$ generated by the ball $B(x_k,k)$, this implies that
\begin{align}\label{L1}
4D^k\geq 1+\sum\limits_{i=0}^{k-1}D(D-1)^i\geq|V|
\end{align}

If $L=2k-1$ for some positive integer $k$, then $d(x_0,x_k)=k>d(x_k,x_L)=k-1$. One easily shows that $|V|$ is less than the subgraph in $T_D$ generated by the ball $B(x_k,k)$, this gives that
\begin{align}\label{L2}
4D^k\geq 1+\sum\limits_{i=0}^{k-1}D(D-1)^i\geq|V|
\end{align}

It follows that $L\geq 2\log_{D}\dfrac{|V|}{4}$ from (\ref{L1}), (\ref{L2}).

\pfd

The following is a corollary of above results.
\pf[Proof of Corollary~\ref{2tend0}]
It follows from Theorem \ref{2d} and Proposition \ref{L}.

\pfd

\rmk For the  subsets $\Omega_n$ in a homogeneous tree, i.e. the tree of constant degree, satisfying with $|\Omega_n|\to\infty,$ $n\to\infty,$ the Steklov problem on $\Omega_n$ satisfying
$$\lambda_2(\Omega_n)\to 0, n\to\infty.$$ A similar result was proved for subgraphs in integer lattices by \cite{Han19}, which was extended to Cayley graphs of polynomial growth by \cite{Perrin20}. The above result is a generalization for Cayley graphs of free groups.
\rmkd

Note that there are some results on lower bounds for general graphs of bounded degree; see \cite{Perrin19}.

In the following, we will investigate the upper bounds of $\lambda_k(\g),$ $k\geq 3,$ for finite trees of bounded degree $D.$ %On account of the multiplicity of eigenvalues, the lower bounds of $\lambda_k(\g)$ are not interesting. We mainly focus on the upper bounds of $\lambda_k(\g)$. The precious statement is as follows.
First, we propose a crucial lemma on the combinatorial property of any finite tree of bounded degree.
\begin{lemma}\label{ktree}
Given a finite tree $\g=(V,E)$ of bounded degree $D$ with boundary $\dO$. Then there admit $k-1$ disjoint subtrees $\g_1,\g_2,\cdots,\g_{k-1}$ in $\g$ such that
$$\dfrac{1}{(D-1)(k-1)}\leq \dfrac{|\dO_j|}{|\dO|}\leq \dfrac{1}{k-1}$$
for $1\leq j\leq k-1$, where $\dO_j$ is the relative boundary of $\g_j$ in $\g$.

\end{lemma}
\pf
Similar to the proof of Lemma~\ref{2tree}, one can take arbitrary edge $e_1\in E$ and remove it to get two connected components of $\g$. Then at least one of the connected components $\h_1$ with the relative boundary $\dO(\h_1)$ satisfying $$\dfrac{|\dO(\h_1)|}{|\dO|}\geq \dfrac{1}{k-1}.$$ There are at most $D-1$ connected components if one removes edges in $\h_1$ adjacent to $e_1$.

We may assume $\h_2$ is the connected component which is obtained by removing $e_2$, with maximal cardinality of the relative boundary $\dO(\h_2)$ among these connected components. If $\dfrac{|\dO(\h_2)|}{|\dO|}<\dfrac{1}{k-1}$, we have $$\dfrac{1}{(D-1)(k-1)}\leq\dfrac{|\dO(\h_1)|}{(D-1)|\dO|}\leq\dfrac{|\dO(\h_2)|}{|\dO|}<\dfrac{1}{k-1}.$$

Otherwise, one can get $$\dfrac{|\dO(\h_2)|}{|\dO|}\geq \dfrac{1}{k-1}.$$ There are at most $D-1$ connected components if one removes edges in $\h_2$ adjacent to $e_2$ and repeat the above procedure. Since $\g$ is finite, we finally get finitely many subgraphs $\h_1\supset\h_2\supset\cdots\supset\h_s$ and consecutive edges $e_1,e_2,\cdots,e_s$, such that $\h_{i}$ is the connected component obtained by removing $e_i$ in $\h_{i-1}$ and $$\dfrac{|\dO(\h_i)|}{|\dO|}\geq \dfrac{1}{k-1}$$ for $\ 1\leq i\leq s-1,$ and $$\dfrac{1}{(D-1)(k-1)}\leq\dfrac{|\dO(\h_s)|}{|\dO|}\leq \dfrac{1}{k-1}.$$

Thus we get a finite,connected subtree $\g_1:=\h_s$ of two connected components by removing $\mu_1:=e_s$ in $\g$ with $$\dfrac{1}{(D-1)(k-1)}\leq\dfrac{|\dO_1|}{|\dO|}=\dfrac{|\dO(\h_s)|}{|\dO|}\leq \dfrac{1}{k-1}.$$ For another connected component $\hat\g_1$ with bigger cardinality of the relative boundary, we use the above method to get a finite, connected subtree $\g_2$ of two connected components by removing $\mu_2$ in $\hat\g_1$ with $$\dfrac{1}{(D-1)(k-1)}\leq\dfrac{|\dO_2|}{|\dO|}\leq \dfrac{1}{k-1}.$$

Repeat this procedure, we can get $k-1$ finite, connected disjoint subtrees $\g_1,\g_2,\cdots,\g_{k-1}$ which are obtained by removing $k-1$ distinct edges $\mu_1,\mu_2,\cdots,\mu_{k-1}$, such that $$\dfrac{1}{(D-1)(k-1)}\leq \dfrac{|\dO_j|}{|\dO|}\leq \dfrac{1}{k-1}$$ for $1\leq j\leq k-1$.

\pfd

First, we prove one of the main results, Theorem \ref{ulamk}.
\pf[Proof of Theorem~\ref{ulamk}]

First, we claim that one can construct $k-1$ functions $f_1,f_2,\cdots,f_{k-1}\in \rv$ with $supp(f_i)\cap supp(f_j)=\emptyset$, $supp(\nabla f_i)\cap supp(\nabla f_j)=\emptyset$ for any $1\leq i\neq j\leq k-1$. As a consequence, $R(\sum\limits_{j=1}^{k-1}b_jf_j)\leq \max\limits_{1\leq j\leq k-1}R(f_j)$ for any $b_1,b_2,\cdots,b_{k-1}\in\R$ with $0\neq f=\sum\limits_{j=1}^{k-1}b_jf_j$.

For the subgraph $\g_j$,
%apply the above result for $k=3$ to two subgraphs $\g_{j,1}=(V_{j,1},E_{j,1})\subset \g_j=(V_j,E_j)$ and $\g_{j,2}=(V_{j,1},E_{j,1}):=\g_j-\g_{j,1}$ with relative boundaries $\dO_{j,1},\dO_{j,2}$ in $\g_j$, such that
repeating the process as
before in the proof of Lemma \ref{ktree} for $k=3$ to two subgraphs $\g_{j,1}=(V_{j,1},E_{j,1})\subset \g_j=(V_j,E_j)$ and $\g_{j,2}=(V_{j,1},E_{j,1}):=\g_j-\g_{j,1}$ with relative boundaries $\dO_{j,1},\dO_{j,2}$ in $\g_j$, one can obtain
\begin{align}\label{do12}
\dfrac{1}{2(D-1)}\leq \dfrac{|\dO_{j,1}|}{|\dO_j|}\leq \dfrac{1}{2}
\leq \dfrac{|\dO_{j,2}|}{|\dO_j|}=1-\dfrac{|\dO_{j,1}|}{|\dO_j|}\leq 1-\dfrac{1}{2(D-1)},
\end{align}
where $1\leq i\leq 2,\ 1\leq j\leq k-1$.

Denote by $\mu_{j,1}\in E_j$ the edge, such that the subgraphs $\g_{j,1},\g_{j,2}$ are obtained from two connected components by removing $\mu_{j,1}$ in $\g_j$.

We define $f_j\in\rv$ in the following.
\begin{equation}\label{f_j}
\left\{
\begin{array}{lr}
f_j(x)=\dfrac{|\dO_{j,2}|}{|\dO_j|},\ x\in V_{j,1}
&\\f_j(x)=-\dfrac{|\dO_{j,1}|}{|\dO_j|},\ x\in V_{j,2}
&\\f_j(x)=0,\ \ elsewhere.
\end{array}
\right.
\end{equation}

By the construction, One easily sees that $(f_j,1_{\dO})=\sum\limits_{x\in\dO}f(x)=0$, $supp(f_j)\subset V_j$ and $supp(\nabla f_j)\subset \{\mu_{j,1},\mu_j\}$ for $1\leq j\leq k-1$. Hence we get the above claim.

Let $W=Span\{f_1,f_2,\cdots,f_{k-1}\}\subset \rv$, by the above claim we have $dim W=k-1$ and $(f,1_{\dO})=0$ for any $0\neq f\in W$. We may assume the previous $f=\sum\limits_{j=1}^{k-1}a_jf_j$. According to (\ref{var k}) and the above claim, we deduce
\begin{align}\label{ulk}
\lambda_2(\g)\leq \max\limits_{0\neq f\in W}R(f)=\max\limits_{1\leq j\leq k-1}R(f_j)
\end{align}

On the other hand, using Lemma \ref{ktree}, the above claim, (\ref{do12}), (\ref{f_j}), we have
\begin{align*}
R(f_j)&=\dfrac{\sum\limits_{(x,y)\in supp(\nabla f_j)}(f_j(x)-f_j(y))^2}{\sum\limits_{x\in \delta\Omega}f_j^2(x)}
\\&\leq \dfrac{\left(\dfrac{|\dO_{j,2}|}{|\dO_j|}\right)^2+\left(\dfrac{|\dO_{j,1}|}{|\dO_j|}+\dfrac{|\dO_{j,2}|}{|\dO_j|}\right)^2}{|\dO_{j,1}|\left(\dfrac{|\dO_{j,2}|}{|\dO_j|}\right)^2+|\dO_{j,2}|\left(\dfrac{|\dO_{j,1}|}
{|\dO_j}|\right)^2}
\\&\leq \dfrac{2}{|\dO|\dfrac{|\dO_{j,1}|}{|\dO_j|}\dfrac{|\dO_{j,2}|}{|\dO_j|}\dfrac{|\dO_{j}|}{|\dO|}}
\\&\leq \dfrac{8(D-1)^2(k-1)}{|\dO|}.
\end{align*}

\pfd

%\rmk
%Note that the constant term of upper bound in Theorem~\ref{ulamk} depends only on $k,D$ and is linear in $k$.
%\rmkd

At the end, we prove the upper bound of $k$-th Steklov eigenvalue using the size of the set of vertices for trees of degree at least three.

\pf[Proof of Corollary \ref{kV}]
It follows from Lemma \ref{dv} and Theorem \ref{ulamk}.
\pfd

%\tm[Theorem~{1.2} in \cite{Hua20}]
%\tmd

%This yields the corollary of Theorem~\ref{thm:main1}.

%\co Let $G$ be an infinite penny graph with bounded facial degree.
%Then for any $k\geq 1,$$$\dim\mathcal{P}^{k}(G)\leq C(D) k^2.$$
%\cod
%\begin{remark} If Conjecture~\ref{con:de} holds, then the above result will be improved to that for any $k\geq 1,$$$\dim\mathcal{P}^{2k}(G)\leq C(D) k^2.$$
%\end{remark}
%This yields the following corollary.
%\begin{corollary}[\cite{Huaancient19}] If $G$ is one of the following:
%\begin{enumerate}
%\item a planar graph with nonnegative combinatorial curvature,
%\item a Cayley graph of polynomial volume growth, or
%\item a graph satisfying $CDE'(0,n),$ $n<\infty,$ \eqref{def:assw} and \eqref{def:delta},
%\end{enumerate}
%then $\dim\mathcal{P}^{k}(G)<\infty.$
%\end{corollary}

%%%%%%%%%%%%%%%%%%%%%%%%%%%%%
%%%%%%%%%%%%%%%%%%%%%%%%%%%%%
%%%%%%%%%%%%%%%%%%%%%%%%%%%%%

\bibliography{py1(1)}
\bibliographystyle{alpha}

\end{document}